# Fibonacci, Lucas, and Spread Polynomials

Johann Cigler

**Abstract.**

This note gives an elementary exposition of a variant of the spread polynomials in terms of Fibonacci and Lucas polynomials.

## 1. Introduction

Norman Wildberger's spread polynomials $S_n(x)$ can be defined as

$$(1) \qquad S_n(x) = \frac{1 - T_n(1 - 2x)}{2},$$

where $T_n(x)$ are the Chebyshev polynomials of the first kind. They can be characterized by the property

$$(2) \qquad S_n(\sin^2 \theta) = \sin^2(n\theta).$$

Hans-Christian Herbig [3] studied the variant

$$(3) \qquad Z_n(x) = 4 S_n\left(\frac{x}{4}\right),$$

which satisfies

$$(4) \qquad Z_n(4\sin^2 \theta) = 4\sin^2(n\theta).$$

Herbig's variant can be best described in terms of different versions of the Fibonacci and Lucas polynomials. We first recall some well-known facts about these polynomials which later will be needed.

We shall write $F_n(x)$ for the version of Fibonacci polynomials defined by the recursion

$$(5) \qquad F_n(x) = x F_{n-1}(x) + F_{n-2}(x)$$

with initial values $F_0(x) = 0$ and $F_1(x) = 1$.

If

$$(6) \qquad \alpha = \alpha(x) = \frac{x + \sqrt{x^2 + 4}}{2}, \ \bar{\alpha} = \bar{\alpha}(x) = \frac{x - \sqrt{x^2 + 4}}{2}$$

denote the roots of the characteristic equation $z^2 - xz - 1 = 0$ of the recursion, we get Binet's formula

(7)
$$F_n(x) = \frac{\alpha^n - \bar{\alpha}^n}{\alpha - \bar{\alpha}}.$$

Note that $\alpha - \bar{\alpha} = \sqrt{x^2 + 4}$ and $\alpha\bar{\alpha} = -1$.

The first terms are

$$(F_n(x))_{n \geq 0} = (0, 1, x, 1 + x^2, 2x + x^3, 1 + 3x^2 + x^4, 3x + 4x^3 + x^5, \cdots).$$

It should be noted that $F_n(x)$ is a monic polynomial of degree $\deg F_n(x) = n - 1$.

For $x = 1$ we get the Fibonacci numbers $(F_n)_{n \geq 0} = (0, 1, 1, 2, 3, 5, 8, 13, \cdots)$.

Let us also mention that the identities

$$\begin{pmatrix} F_{n-1}(x) & F_n(x) \\ F_n(x) & F_{n+1}(x) \end{pmatrix} \begin{pmatrix} 0 & 1 \\ 1 & x \end{pmatrix} = \begin{pmatrix} F_n(x) & F_{n+1}(x) \\ F_{n+1}(x) & F_{n+2}(x) \end{pmatrix} \text{ and } \begin{pmatrix} F_0(x) & F_1(x) \\ F_1(x) & F_2(x) \end{pmatrix} = \begin{pmatrix} 0 & 1 \\ 1 & x \end{pmatrix}$$

imply

(8)
$$\begin{pmatrix} F_{n-1}(x) & F_n(x) \\ F_n(x) & F_{n+1}(x) \end{pmatrix} = \begin{pmatrix} 0 & 1 \\ 1 & x \end{pmatrix}^n.$$

Applying determinants, this gives Cassini's formula

(9)
$$F_{n-1}(x)F_{n+1}(x) - F_n^2(x) = (-1)^n.$$

Similarly, we write $L_n(x)$ for the version of Lucas polynomials which are defined by

(10)
$$L_n(x) = xL_{n-1}(x) + L_{n-2}(x)$$

with $L_0(x) = 2$ and $L_1(x) = x$.

Here we get

(11)
$$L_n(x) = \alpha^n + \bar{\alpha}^n.$$

The first terms are

$$(L_n(x))_{n \geq 0} = (2, x, 2 + x^2, 3x + x^3, 2 + 4x^2 + x^4, 5x + 5x^3 + x^5, 2 + 9x^2 + 6x^4 + x^6, \cdots).$$

Note that $L_n(x)$ is a monic polynomial of degree $\deg L_n(x) = n$.

For $x = 1$ we get the Lucas numbers

$$(L_n)_{n \geq 0} = (2, 1, 3, 4, 7, 11, 18, \cdots).$$

We shall write $f_n(x)$ for the version of Fibonacci polynomials defined by

(12) $$f_n(x) = x f_{n-1}(x) - f_{n-2}(x)$$

with $f_0(x) = 0$ and $f_1(x) = 1$.

Let

(13) $$\beta = \beta(x) = \frac{x + \sqrt{x^2 - 4}}{2},$$

$$\bar{\beta} = \bar{\beta}(x) = \frac{x - \sqrt{x^2 - 4}}{2}$$

be the roots of the characteristic equation $z^2 - xz + 1 = 0$.

Then we get

(14) $$f_n(x) = \frac{\beta^n - \bar{\beta}^n}{\beta - \bar{\beta}}.$$

Note that $\beta - \bar{\beta} = \sqrt{x^2 - 4}$ and $\beta \bar{\beta} = 1$.

The first terms are

$$(f_n(x))_{n \geq 0} = (0, 1, x, -1 + x^2, -2x + x^3, 1 - 3x^2 + x^4, 3x - 4x^3 + x^5, -1 + 6x^2 - 5x^4 + x^6, \cdots).$$

Note that $f_n(x)$ is a monic polynomial of degree $\deg f_n(x) = n - 1$.

Similarly, we write $l_n(x)$ for the version of Lucas polynomials defined by

(15) $$l_n(x) = x l_{n-1}(x) - l_{n-2}(x)$$

with $l_0(x) = 2$ and $l_1(x) = x$.

Here we get

(16) $$l_n(x) = \beta^n + \bar{\beta}^n.$$

The first terms are

$$(l_n(x))_{n \geq 0} = (2, x, -2 + x^2, -3x + x^3, 2 - 4x^2 + x^4, 5x - 5x^3 + x^5, -2 + 9x^2 - 6x^4 + x^6, \cdots).$$

Note that $l_n(x)$ is a monic polynomial of degree $\deg l_n(x) = n$.

**Remark**

The polynomials $f_n(x)$ and $l_n(x)$ are related to the Chebyshev polynomials $U_n(x)$ and $T_n(x)$ by

(17)
$$f_n(x) = U_{n-1}\left(\frac{x}{2}\right),$$
$$l_n(x) = 2T_n\left(\frac{x}{2}\right).$$

Since $\beta_n(2\cos\theta) = e^{i\theta}$ and $\bar{\beta}_n(2\cos\theta) = e^{-i\theta}$ we get

(18)
$$f_n(2\cos\theta) = \frac{\sin(n\theta)}{\sin\theta},$$
$$l_n(2\cos\theta) = 2\cos(n\theta).$$

Note that

(19) $$l_n\left(x + \frac{1}{x}\right) = x^n + \frac{1}{x^n}.$$

This follows from the recursion $x^n + \frac{1}{x^n} = \left(x + \frac{1}{x}\right)\left(x^{n-1} + \frac{1}{x^{n-1}}\right) - \left(x^{n-2} + \frac{1}{x^{n-2}}\right)$ and the initial values $2$ and $x + \frac{1}{x}$.

Since $x^{2n} + \frac{1}{x^{2n}} = \left(x^n + \frac{1}{x^n}\right)^2 - 2$ we get

(20) $$l_{2n}(x) = l_n^2(x) - 2.$$

Another important property is

(21) $$l_m \circ l_n(x) = l_m(l_n(x)) = l_{mn}(x).$$

It suffices to verify that for $x = z + \frac{1}{z}$.

$$l_{mn}\left(z+\frac{1}{z}\right) = z^{mn} + \frac{1}{z^{mn}} = \left(z^n\right)^m + \frac{1}{\left(z^n\right)^m} = l_m\left(z^n + \frac{1}{z^n}\right) = l_m\left(l_n\left(z+\frac{1}{z}\right)\right).$$

From $\beta^n - \bar{\beta}^n = \sqrt{x^2 - 4} f_n(x)$ and $\beta^n + \bar{\beta}^n = l_n(x)$ we get

(22)
$$\beta^n = \frac{l_n(x) + \sqrt{x^2 - 4} f_n(x)}{2},$$
$$\bar{\beta}^n = \frac{l_n(x) - \sqrt{x^2 - 4} f_n(x)}{2},$$
$$l_n^2(x) - \left(x^2 - 4\right) f_n^2(x) = 4\beta\bar{\beta} = 4.$$

For $x = 2\cos\theta$ these identities reduce to $\cos\theta \pm i\sin\theta = e^{\pm i\theta}$ and $\cos^2\theta + \sin^2\theta = 1$.

**2. The polynomials $Z_n(x)$.**

From (1) and (17) we get

(23) $$Z_n(x) = 2 - l_n(2 - x).$$

The first terms are

$$(Z_n(x))_{n\geq 0} = \left(0, x, 4x - x^2, 9x - 6x^2 + x^3, 16x - 20x^2 + 8x^3 - x^4, 25x - 50x^2 + 35x^3 - 10x^4 + x^5, \cdots\right).$$

Note that $(-1)^n Z_n(-x)$ is a monic polynomial of degree $\deg Z_n(x) = n$.

The polynomials $Z_n(x)$ satisfy

(24)
$$Z_{mn}(x) = Z_n(Z_m(x)),$$
$$Z_n\left(-\left(x - \frac{1}{x}\right)^2\right) = -\left(x^n - \frac{1}{x^n}\right)^2.$$

This follows from $Z_{mn}(x) = 2 - l_{nm}(2 - x) = 2 - l_n\left(l_m(2 - x)\right) = 2 - l_n\left(2 - Z_m(x)\right) = Z_n(Z_m(x))$
and

$$Z_n\left(-\left(x - \frac{1}{x}\right)^2\right) = 2 - l_n\left(2 + \left(x - \frac{1}{x}\right)^2\right) = 2 - l_n\left(x^2 + \frac{1}{x^2}\right) = 2 - l_n\left(l_2\left(x + \frac{1}{x}\right)\right)$$

$$= 2 - l_{2n}\left(x + \frac{1}{x}\right) = 2 - \left(x^{2n} + \frac{1}{x^{2n}}\right) = -\left(x^n - \frac{1}{x^n}\right)^2.$$

Observing that

$$l_n(2-x) = \lambda^n + \bar{\lambda}^n$$

(25) with $\lambda = \lambda(x) = \beta(2-x) = \dfrac{2-x+\sqrt{x^2-4x}}{2}$,

$$\bar{\lambda} = \bar{\lambda}(x) = \dfrac{2-x-\sqrt{x^2-4x}}{2},$$

we see that the polynomials $Z_n(x)$ satisfy the linear recursion

(26) $$Z_n(x) = (3-x)Z_{n-1}(x) + (x-3)Z_{n-2}(x) + Z_{n-3}(x),$$

because by (23) they are linear composita of $\lambda^n(x)$, $\bar{\lambda}^n(x)$, and $1^n$ and
$(z-\lambda)(z-\bar{\lambda})(z-1) = (z^2 - (2-x)z + 1)(z-1) = z^3 - (3-x)z^2 + (3-x)z - 1.$

For $x \in \{0,1,2,3,4\}$ $\lambda(x)$ is a primitive root of unity. Therefore, the sequence $(Z_n(x))$ is periodic for these values.

Let

(27) $$\mu = \dfrac{\sqrt{x-4}+\sqrt{x}}{2} = \alpha(\sqrt{x-4}),$$

$$\bar{\mu} = \bar{\alpha}(\sqrt{x-4}) = \dfrac{\sqrt{x-4}-\sqrt{x}}{2}.$$

Then

(28) $$\mu^2 = -\bar{\lambda},$$
$$\bar{\mu}^2 = -\lambda,$$
$$\mu\bar{\mu} = -1.$$

This implies

**Theorem 1**

*The polynomials $Z_n(x)$ satisfy*

(29) $$Z_n(x) = (-1)^{n-1} x \left(F_n\left(\sqrt{x-4}\right)\right)^2$$

*and*

(30) $$Z_n(x) = -\dfrac{(\lambda(x)^n - 1)^2}{\lambda(x)^n}.$$

**Proof**

Identity (29) follows from

$$Z_n(x) = 2 - l_n(2-x) = 2 - \left(\lambda^n + \bar{\lambda}^n\right) = 2 - \left((-\bar{\mu}^2)^n + (-\mu^2)^n\right) = 2 - (-1)^n\left(\mu^{2n} + \bar{\mu}^{2n}\right)$$

$$= (-1)^{n-1}\left(\mu^n - \bar{\mu}^n\right)^2 = (-1)^{n-1} x \left(\frac{\mu^n - \bar{\mu}^n}{\mu - \bar{\mu}}\right)^2 = (-1)^{n-1} x \left(F_n\left(\sqrt{x-4}\right)\right)^2.$$

For $x = 5$ we get squares of Fibonacci numbers

$$\left(Z_n(5)\right)_{n\geq 0} = 5(0, 1, -1, 4, -9, 25, -64, 169, -441, \cdots) = 5\left((-1)^{n-1} F_n^2\right)_{n\geq 0}.$$

Identity (30) follows from

$$Z_n(x) = 2 - l_n(2-x) = 2 - \left(\lambda^n + \bar{\lambda}^n\right) = 2 - \lambda^n - \frac{1}{\lambda^n} = -\frac{\left(\lambda^n - 1\right)^2}{\lambda^n} = -\frac{\left(\bar{\lambda}^n - 1\right)^2}{\bar{\lambda}^n}$$

**Proposition 1**

(31)
$$F_{2n+1}\left(\sqrt{x-4}\right) = \frac{l_{2n+1}\left(\sqrt{x}\right)}{\sqrt{x}},$$

$$F_{2n}\left(\sqrt{x-4}\right) = \frac{f_{2n}\left(\sqrt{x}\right)}{\sqrt{x}} \sqrt{x-4}.$$

**Proof**

$f_n\left(\sqrt{x}\right)$ and $l_n\left(\sqrt{x}\right)$ are linear composita of $\beta^n\left(\sqrt{x}\right)$ and $\bar{\beta}^n\left(\sqrt{x}\right)$ and $F_n\left(\sqrt{x-4}\right)$ is a linear compositum of $\mu^n$ and $\bar{\mu}^n$.

Observing that $\beta\left(\sqrt{x}\right) = \frac{\sqrt{x} + \sqrt{x-4}}{2} = \mu(x)$ and $\bar{\beta}\left(\sqrt{x}\right) = \frac{\sqrt{x} - \sqrt{x-4}}{2} = -\bar{\mu}(x)$ we get

$$F_{2n+1}\left(\sqrt{x-4}\right) = \frac{\mu^{2n+1}(x) - \bar{\mu}^{2n+1}(x)}{\mu(x) - \bar{\mu}(x)} = \frac{\beta^{2n+1}\left(\sqrt{x}\right) + \bar{\beta}^{2n+1}\left(\sqrt{x}\right)}{\sqrt{x}} = \frac{l_{2n+1}\left(\sqrt{x}\right)}{\sqrt{x}}$$

and

$$F_{2n}\left(\sqrt{x-4}\right) = \frac{\mu^{2n}(x) - \bar{\mu}^{2n}(x)}{\mu(x) - \bar{\mu}(x)} = \frac{\beta^{2n}\left(\sqrt{x}\right) - \bar{\beta}^{2n}\left(\sqrt{x}\right)}{\sqrt{x}}$$

$$= \frac{\beta^{2n}\left(\sqrt{x}\right) - \bar{\beta}^{2n}\left(\sqrt{x}\right)}{\beta\left(\sqrt{x}\right) - \bar{\beta}\left(\sqrt{x}\right)} \frac{\beta\left(\sqrt{x}\right) - \bar{\beta}\left(\sqrt{x}\right)}{\sqrt{x}} = f_{2n}\left(\sqrt{x}\right)\frac{\sqrt{x-4}}{\sqrt{x}}.$$

This implies

(32)
$$Z_{2n+1}(x) = xF_{2n+1}^2\left(\sqrt{x-4}\right) = l_{2n+1}^2\left(\sqrt{x}\right),$$
$$Z_{2n}(x) = -xF_{2n}^2\left(\sqrt{x-4}\right) = f_{2n}^2\left(\sqrt{x}\right)(4-x).$$

The polynomials $Z_n(x)$ also satisfy

(33) $$Z_{n-1}(x)Z_{n+1}(x) = (Z_n(x) - x)^2$$

because by Cassini's formula (9)

$$Z_{n-1}(x)Z_{n+1}(x) = x^2\left(F_{n-1}^2\left(\sqrt{x-4}\right)F_{n+1}^2\left(\sqrt{x-4}\right)\right) = x^2\left(F_n^2\left(\sqrt{x-4}\right) + (-1)^n\right)^2$$
$$= \left(xF_n^2\left(\sqrt{x-4}\right) + (-1)^n x\right)^2 = \left((-1)^{n-1}Z_n(x) - (-1)^{n-1}x\right)^2 = (Z_n(x) - x)^2.$$

The following Theorem is a refinement of a conjecture of Shuxiang Goh and N. J. Wildberger which has been mentioned in [3] and first proved in [2] with other methods.

**Theorem 2**

*There is a sequence of monic polynomials $(\Phi_n(x))_{n\geq 1}$, $\Phi_n(x) \in \mathbb{Z}[x]$, such that $(-1)^{n-1}Z_n(x) = \prod_{d|n}\Phi_d(x)$ and $\deg \Phi_d(x) = \varphi(d)$, where $\varphi(d)$ is Euler's totient function*

*If $\phi_n(x)$ denotes the minimal polynomial of $4\sin^2\left(\frac{\pi}{n}\right)$ then $\Phi_1(x) = \phi_1(x) = x$, $\Phi_2(x) = \phi_2(x) = x - 4$, and $\Phi_n(x) = \phi_n(x)^2$ for $n \geq 3$.*

**Proof**

Let

(34) $$C_n(x) = \prod_{\substack{1 \leq k \leq n \\ \gcd(k,n)=1}} \left(x - e^{2\pi i \frac{k}{n}}\right)$$

be the cyclotomic polynomials. They satisfy

(35) $$x^n - 1 = \prod_{d|n} C_d(x).$$

The first cyclotomic polynomials are

$$(C_n(x))_{n\geq 1} = (-1+x, 1+x, 1+x+x^2, 1+x^2, 1+x+x^2+x^3+x^4, 1-x+x^2, \cdots)$$

By (30) $Z_n(x) = -\dfrac{(\lambda(x)^n - 1)^2}{\lambda(x)^n}$. This together with (35) and the identity $\sum_{d|n} \varphi(d) = n$ for Euler's totient function $\varphi(n)$ gives

(36) $$-Z_n(x) = \dfrac{\prod_{d|n} C_d(\lambda(x))^2}{\lambda(x)^n} = \prod_{d|n} \dfrac{C_d(\lambda(x))^2}{\lambda(x)^{\varphi(d)}}.$$

Let us simplify the right-hand side.

$$\dfrac{C_1(\lambda(x))^2}{\lambda(x)^{\varphi(1)}} = \dfrac{(\lambda(x)-1)^2}{\lambda(x)} = \lambda(x) + \dfrac{1}{\lambda(x)} - 2 = \lambda(x) + \bar{\lambda}(x) - 2 = 2 - x - 2 = -x,$$

$$\dfrac{C_2(\lambda(x))^2}{\lambda(x)^{\varphi(2)}} = \dfrac{(\lambda(x)+1)^2}{\lambda(x)} = \lambda(x) + \bar{\lambda}(x) + 2 = 4 - x.$$

For $n \geq 3$ the cyclotomic polynomial $C_n(x)$ is a monic symmetric polynomial with integer coefficients of even degree $\varphi(n)$. Therefore $\dfrac{C_n(x)}{x^{\frac{\varphi(n)}{2}}} = \sum_{k=0}^{\frac{\varphi(n)}{2}} c_k \left(x^k + \dfrac{1}{x^k}\right)$ is a sum of terms $x^k + \dfrac{1}{x^k}$ with integer coefficients $c_k$. Since $\lambda(x)^k + \dfrac{1}{\lambda(x)^k} = l_k(2-x)$ we get

$$\psi_n(x) := (-1)^{\frac{\varphi(n)}{2}} \dfrac{C_n(\lambda(x))}{(\lambda(x))^{\frac{\varphi(n)}{2}}} = (-1)^{\frac{\varphi(n)}{2}} \sum_{k=0}^{\frac{\varphi(n)}{2}} c_k l_k(2-x) \text{ is a monic polynomial of degree}$$

$\deg \psi_n(x) = \dfrac{\varphi(n)}{2}$ with integer coefficients $c_k$.

For example,

$C_3(x) = 1 + x + x^2$ gives $\dfrac{C_3(x)}{x} = 1 + \left(x + \dfrac{1}{x}\right)$ and $\psi_3(x) = -(1 + l_1(2-x)) = -(3-x) = x - 3,$

$C_4(x) = 1 + x^2$ gives $\dfrac{C_4(x)}{x} = x + \dfrac{1}{x}$ and $\psi_4(x) = -l_1(2-x) = x - 2,$

$C_5(x) = 1 + x + x^2 + x^3 + x^4$ gives $\dfrac{C_5(x)}{x^2} = 1 + \left(x + \dfrac{1}{x}\right) + \left(x^2 + \dfrac{1}{x^2}\right)$ and

$\psi_5(x) = 1 + l_1(2-x) + l_2(2-x) = 1 + (2-x) + (2 - 4x + x^2) = x^2 - 5x + 5.$

Since

$$\lambda(4\sin^2\theta) = \dfrac{2 - 4\sin^2\theta + 4\sin\theta\sqrt{(\sin^2\theta - 1)}}{2} = \dfrac{2\cos(2\theta) + 4i\sin\theta\cos\theta}{2} = e^{2i\theta}$$

we get $C_n\left(\lambda\left(4\sin^2\dfrac{\pi}{n}\right)\right) = C_n\left(e^{\frac{2\pi i}{n}}\right) = 0.$

Thus $\psi_n(x) = (-1)^{\frac{\varphi(n)}{2}} \dfrac{C_n(\lambda(x))}{\lambda(x)^{\frac{\varphi(n)}{2}}}$ is a monic polynomial with integer coefficients of degree $\dfrac{\varphi(n)}{2}$ with root $4\sin^2\left(\dfrac{\pi}{n}\right)$. Therefore it coincides with the minimal polynomial $\phi_n(x)$ of $4\sin^2\left(\dfrac{\pi}{n}\right)$.

From
$(\phi_n(x))_{n \geq 1} = \left(x, -4 + x, -3 + x, -2 + x, 5 - 5x + x^2, -1 + x, -7 + 14x - 7x^2 + x^3, 2 - 4x + x^2, \cdots\right)$ we see that

(37)
$$\begin{aligned}\Phi_1(x) &= x = \phi_1(x), \\ \Phi_2(x) &= x - 4 = \phi_2(x), \\ \Phi_n(x) &= \phi_n^2(x) \text{ for } n \geq 3.\end{aligned}$$

For example,

$Z_1(x) = x = \Phi_1(x)$,  $-Z_2(x) = x^2 - 4x = x(x-4) = \Phi_1(x)\Phi_2(x),$

$Z_3(x) = x^3 - 6x^2 + 9x = x(x-3)^2 = \Phi_1(x)\Phi_3(x),$

$-Z_4(x) = x^4 - 8x^3 + 20x^2 - 16x = x(x-4)(x-2)^2 = \Phi_1(x)\Phi_2(x)\Phi_4(x).$

**Theorem 3**

*For odd m we get*

(38) $$Z_{2m}(x) = Z_m(x) Z_m(-\Phi_2(x))$$

*and for* $k \geq 2$

(39) $$Z_{2^k m}(x) = Z_{2^{k-1} m}(x) Z_m\left(\Phi_{2^k}(x)\right).$$

**Proof**

$$Z_{2m}(x) = -\frac{\left(\lambda(x)^{2m}-1\right)^2}{\lambda(x)^{2m}} = -\frac{\left(\lambda(x)^m-1\right)^2}{\lambda(x)^m} \frac{\left(\lambda(x)^m+1\right)^2}{\lambda(x)^m} = Z_m(x) \frac{\left(\lambda(x)^m+1\right)^2}{\lambda(x)^m}.$$

Since $\lambda(4-x) = \dfrac{2-(4-x)+\sqrt{(4-x)^2-4(4-x)}}{2} = \dfrac{x-2+\sqrt{x^2-4}}{2} = -\bar{\lambda}(x)$

We get

$$Z_m(-\Phi_2(x)) = Z_m(4-x) = -\frac{\left(\lambda(4-x)^m-1\right)^2}{\lambda(4-x)^m} = \frac{\left(\bar{\lambda}(x)^m+1\right)^2}{\bar{\lambda}(x)^m}$$

$$= -\bar{\lambda}(x)^m - \lambda(x)^m - 2 = \frac{\left(\lambda(x)^m+1\right)^2}{\lambda(x)^m},$$

which proves (38).

Observing that $Z_{2^k}(x) = \prod_{j=0}^{k} \Phi_{2^j}(x)$ and $\Phi_{2^k}(x) = \phi_{2^k}^2(x)$ for $k \geq 2$ we get

$$\Phi_{2^k}(x) = \frac{Z_{2^k}(x)}{Z_{2^{k-1}}(x)} = \frac{1}{\lambda(x)^{2^{k-1}}} \left(\frac{\lambda^{2^k}(x)-1}{\lambda^{2^{k-1}}(x)-1}\right)^2 = \frac{\left(\lambda^{2^{k-1}}(x)+1\right)^2}{\lambda(x)^{2^{k-1}}} = \left(\lambda^{2^{k-2}}(x) + \frac{1}{\lambda^{2^{k-2}}(x)}\right)^2$$ which

gives

(40) $$\phi_{2^k}(x) = l_{2^{k-2}}\left(\lambda(x) + \frac{1}{\lambda(x)}\right) = l_{2^{k-2}}(2-x)$$

for $k \geq 3$.

From $l_{2^{n+1}}(x) = l_2(l_{2^n}(x)) = l_{2^n}^2(x) - 2$ we get for $k \geq 2$ the recursion $\phi_{2^{k+1}}(x) = \phi_{2^k}^2(x) - 2$ and therefore

(41) $$\Phi_{2^{k+1}}(x) = \left(\Phi_{2^k}(x) - 2\right)^2.$$

Observing that

(42) $$\lambda\left((x-2)^2\right) = -\lambda(x)^2$$

and that $(x-2)^2 = \Phi_4(x)$ we get with induction

(43) $$\lambda\left(\Phi_{2^k}(x)\right)=-\lambda(x)^{2^{k-1}}.$$

This gives

$$Z_{2^k m}(x) = -\frac{\left(\lambda(x)^{2^k m}-1\right)^2}{\lambda(x)^{2^k m}} = -\frac{\left(\lambda(x)^{2^{k-1} m}-1\right)^2}{\lambda(x)^{2^{k-1} m}} \frac{\left(\lambda(x)^{2^{k-1} m}+1\right)^2}{\lambda(x)^{2^{k-1} m}}$$

$$= -Z_{2^{k-1} m}(x) \frac{\left(\lambda\left(\Phi_{2^k}(x)\right)^m - 1\right)^2}{\lambda\left(\Phi_{2^k}(x)\right)^m} = Z_{2^{k-1} m}(x) Z_m\left(\Phi_{2^k}(x)\right).$$

Finally let us state a product representation of the Fibonacci polynomials.

By (29) $Z_n(x) = x F_n^2\left(\sqrt{x-4}\right)$. Therefore, $F_n^2\left(\sqrt{x-4}\right) = \prod_{\substack{d|n \\ d>1}} \Phi_d(x)$ and

$F_n^2(x) = \prod_{\substack{d|n \\ d>1}} \Phi_d(x^2+4)$. Since $\Phi_2(x^2+4) = x^2$ all terms are squares.

If we define

$\rho_n(x) = \phi_n(x^2+4)$ for $n \geq 3$ and $\rho_2(x) = x$ we get

$$\left(\rho_n(x)\right)_{n \geq 2} = \left(x, 1+x^2, 2+x^2, 1+3x^2+x^4, 3+x^2, 1+6x^2+5x^4+x^6, 2+4x^2+x^4, \cdots\right)$$

and

(44) $$F_n(x) = \prod_{\substack{d|n \\ d>1}} \rho_d(x).$$

**Remark**

It is well known that the Fibonacci numbers and polynomials are strong divisibility sequences.

(An overview about this concept can be found in [1].) Theorem 2 says that the spread polynomials also form a strong divisibility sequence.